\documentclass[reqno]{amsart}
\usepackage{amsmath,amssymb,amsthm}
\usepackage{enumerate}
\usepackage{bbm}
\usepackage[utf8]{inputenc}
\usepackage{xcolor}
\newtheorem{theorem}{Theorem}[section]
\newtheorem{corollary}[theorem]{Corollary}
\newtheorem{proposition}[theorem]{Proposition}
\newtheorem{lemma}[theorem]{Lemma}

\newtheorem*{theoremA}{Theorem A}
\newtheorem*{theoremB}{Theorem B}

\theoremstyle{definition}  
\newtheorem{remark}[theorem]{Remark}
\newtheorem{definition}[theorem]{Definition}

\newtheorem{example}[theorem]{Example}

\newcommand{\D}{\mathbb{D}}
\newcommand{\disk}{\mathbb{D}}
\newcommand{\C}{\mathbb{C}}
\newcommand{\what}[1]{\widehat{\,#1\,}}
\newcommand{\inner}[2]{\langle#1,#2\rangle}
\newcommand{\bfQ}{\mathbf{Q}}
\newcommand{\bfP}{\mathbf{P}}
\newcommand{\bfS}{\mathbf{S}}

\newcommand{\SQH}{\text{SQH\,}}

\DeclareMathOperator{\rank}{rank}

\newcommand{\ds}{\displaystyle}

\numberwithin{equation}{section}

\begin{document}

\title[Finite rank perturbations of Toeplitz products]{Finite Rank Perturbations of Toeplitz Products on the Bergman Space}
 
\author{Trieu Le}
\address{Department of Mathematics and Statistics, The University of Toledo, Toledo, OH 43606}
\email{Trieu.Le2@utoledo.edu}

\author{Damith Thilakarathna}
\address{Department of Mathematics and Statistics, The University of Toledo, Toledo, OH 43606}
\email{Damith.Thilakarathna@utoledo.edu}



\begin{abstract}
In this paper we investigate when a finite sum of products of two Toeplitz operators with quasihomogeneous symbols is a finite rank perturbation of another Toeplitz operator on the Bergman space. We discover a noncommutative convolution $\diamond$ on the space of quasihomogeneous functions and use it in solving the problem. Our main results show that if $F_j, G_j$ ($1\leq j\leq N$) are polynomials of $z$ and $\bar{z}$ then $\sum_{j=1}^{N}T_{F_j}T_{G_j}-T_{H}$ is a finite rank operator for some $L^{1}$-function $H$ if and only if $\sum_{j=1}^{N}F_j\diamond G_j$ belongs to $L^1$ and $H=\sum_{j=1}^{N}F_j\diamond G_j$. In the case $F_j$'s are holomorphic and $G_j$'s are conjugate holomorphic, it is shown that $H$ is a solution to a system of first order partial differential equations with a constraint.
\end{abstract}

\maketitle 
 
\section{Introduction}

Let $\D$ denote the open unit disk on the complex plane and let $dA$ be the normalized Lebesgue area measure on $\D$. The Bergman space $A^2$ consists of all holomorphic functions on $\D$ that are square integrable with respect to $dA$. It is well known that $A^2$ is a closed subspace of $L^2=L^2(\D,dA)$. Let $P$ denote the orthogonal projection from $L^2$ onto $A^2$. For a function  $f \in L^2$, the Toeplitz operator $T_f$ with symbol $f$ is defined as
\[T_{f}(h) = P(fh)\]
for all $h\in A^2$ for which the product $f h$ belongs to $L^2$. Since bounded holomorphic functions belong to its domain, the operator $T_{f}$ is densely defined on $A^2$. 
It is clear that if $f$ is bounded, then $T_{f}$ is bounded and $\|T_{f}\|\leq\|f\|_{\infty}$. On the other hand, unbounded symbols may still give rise to bounded Toeplitz operators. See \cite[Chapter~7]{ZhuAMS2007} for a detailed introduction to Toeplitz operators on the Bergman space.

There is an intensive literature on the study of Toeplitz operators. In this paper we are concerned with the problem of when a product $T_{f}T_{g}$ of two Toeplitz operators  is a finite rank perturbation of another Toeplitz operator. More generally, we are interested in finding conditions on the symbols $f_1,\ldots, f_n$ and $g_1, \ldots, g_n$ for which $T_{f_1}T_{g_1}+\cdots+T_{f_n}T_{g_n}$ is a sum of a finite rank operator and another Toeplitz operator. The problem is wide open for general symbols but we have made some progress for symbols that are finite sums of functions of a special type. In particular, our result applies to the case where the symbols are polynomials in $z$ and $\bar{z}$.

Most of the study of Toeplitz operators has been influenced by the seminal paper \cite{brownhalmos1964} of Brown and Halmos. Among other things, they showed that for Toeplitz operators on the Hardy space over the unit circle, one has $T_{f}T_{g}=T_{h}$ if and only if either $g$ or $\bar{f}$ is holomorphic and $h=fg$. The situation on the Bergman space turns out to be much more complicated. Zheng \cite{ZhengJFA1989} showed that if $f, g$ are bounded harmonic functions such that $T_{f}T_{g}=T_{fg}$ on $A^2$, then either $g$ or $\bar{f}$ is holomorphic. Ahern and {\v{C}}u{\v{c}}kovi{\'c} \cite{aherncuckovic2001} obtained an analogue of Brown--Halmos's result under the assumption that $f, g$ are bounded harmonic and $h$ is a $C^{2}$-function such that $\tilde{\Delta}h=(1-|z|^2)^2\Delta h$ is bounded. Ahern later \cite{ahernJFA2004} removed the conditions on $h$ and established that the conclusion remains true for any bounded function $h$. In \cite{GuoSunZhengIJM2007}, Guo, Sun, and Zheng showed that if the operator $T_{f}T_{g}-T_{fg}$ has finite rank, then the rank must actually be zero. {\v{C}}u{\v{c}}kovi{\'c} \cite{cuckovicIEOT2007} obtained criteria for $T_{f}T_{g}-T_{h^n}$ to have finite rank, where $f,g$ and $h$ are bounded harmonic.  The case of finite sums of products of Toeplitz operators was studied by Choe, Koo and Lee in \cite{ChoeRMI2008} for a single variable and in \cite{ChoePA2009} for several variables over the polydisk. In a recent paper, Ding, Qin and Zheng \cite{ding2017theorem} provided a more complete answer to the possible rank of $T_{f}T_{g}-T_{h}$ under the assumption that $f, g$ are bounded harmonic and $h$ is a $C^2$-function such that $\tilde{\Delta}h$ belongs to $L^1(\D,dA)$. It was showed that if the rank is at most one, then it must actually be zero. On the other hand, for each $m\geq 2$, there are examples of rational functions $f,g,h$ for which the rank equals to $m$. The authors of \cite{ding2017theorem} also provided a complete characterization of the functions $f, g$ and $h$ in this case.

For non-harmonic symbols, Ahern \cite[Corollary 1]{ahernJFA2004} exhibited the following interesting example:
\begin{align}
\label{Ex:Ahern}
T_{z}T_{\bar{z}^2{z}} = T_{2\bar{z}z-1},
\end{align}
which says that the Brown-Halmos's result fails for general Toeplitz operators on the Bergman space. The symbols $f$ and $g$ in \eqref{Ex:Ahern} belong to the class of quasihomogeneous functions. We call a function on the unit disk \textit{quasihomogeneous of degree $m$} if it can be represented in the form $\phi(r)e^{im\theta}$ almost everywhere for some function $\phi$ defined on the interval $[0,1)$.  Louhichi, Zakariasy and Strouse \cite{louhichi2006products} obtained criteria for a  product of two Toeplitz operators with quasihomogeneous symbols to equal another Toeplitz operator.     
In their results, the symbols are allowed to blow up at the origin but they are bounded on a certain annulus $\{z: r<|z|<1\}$. 
\v{C}u\v{c}kovi\'c and Louhichi \cite{cuckovic2008finite} investigated finite rank semicommutators and commutators of Toeplitz operators with quasihomogeneous symbols. They proved that for $m\geq n\geq 0$, if the operator $T_{\phi e^{im\theta}}  T_{\psi e^{-in\theta}}-T_{\phi\psi e^{i(m-n)\theta}}$  has finite rank $N$, then $N$ is at most equal to the quasihomogeneous degree $n$. Example 1 in \cite{cuckovic2008finite} shows that $T_{\bar{z}^{-1}}T_{\bar{z}}-T_{1}$ is a rank one operator. By considering $T_{f}T_{g}-T_{h}$, one obtains many other interesting examples (see \cite[Example 1.2]{ding2017theorem}):
\begin{align}
T_{z^2}T_{\bar{z}^3} - T_{3\bar{z}-\frac{2}{z}} = 1\otimes z,\label{ExampleA}\\
T_{z^2}T_{\bar{z}^3z} - T_{3|z|^2-2} = \frac{1}{2}\otimes 1.\label{ExampleB}
\end{align}
Here for any two non-zero functions $u,v\in A^2$, we use $u\otimes v$ to denote the rank one operator defined by $(u\otimes v)(h) = \inner{h}{v}u$ for $h\in A^2$. As a consequence of our main results, we are able to characterize all such examples. That is, we obtain conditions for which the product $T_{r^pe^{im\theta}}T_{r^qe^{in\theta}}$ is a finite rank perturbation of a Toeplitz operator, where $p, q\geq -1$ are real numbers and $m, n$ are integers. Our result applies to finite sums of such products as well. 

\begin{remark}
Throughout the paper, we encounter finite sums of products of possibly unbounded Toeplitz operators whose domains contain the space of holomorphic polynomials. Such an operator is said to have finite rank if the image of holomorphic polynomials is finite dimensional. Of course if the operator happens to be bounded, then the image of all $A^2$ will be of finite dimensional as well.
\end{remark}

While \cite{ding2017theorem} relied heavily on the Berezin transform, our approach in this paper makes use of the Mellin transform. We define a new kind of convolution $\diamond$ and show that if $f$ and $g$ are combinations of certain quasihomogeneous functions, then $T_{f}T_{g}-T_{h}$ is a finite rank operator if and only if $h=f\diamond g$. We obtain

\begin{theoremA}
\label{T:sum_product_polynomials}
Let  $F_j$ and $G_j$ be polynomials in $z$ and $\bar{z}$ for $j=1,\ldots,N$. Then there exists $H\in L^1(\mathbb{D}, dA)$ such that $\ds \sum_{j=1}^{N}T_{F_j}T_{G_j}-T_{H}$ is of finite rank if and only if $\ds \sum_{j=1}^{N}F_j\diamond G_j\in L^1(\D,dA)$.
\end{theoremA}

It is interesting to note that our convolution is similar to the single variable convolution in the study of Berezin-Toeplitz calculus on the Fock-Segal-Bargmann space discovered by Coburn \cite{CoburnPAMS2001} and later developed to more general settings by Bauer \cite{BauerJFA2009}.

For symbols that are holomorphic polynomials, we prove a surprising and exciting connection of the theory of Toeplitz operators to solutions of first order differential equations. 

\begin{theoremB}
\label{T:TOs_ODEs}
Let  $P_j$ and $Q_j$ be Laurent polynomials for $j=1,\ldots, N$ and let $H$ belong to  $L^1(\D,dA)$. Then $\ds\sum_{j=1}^{N}T_{P_j}T_{\overline{Q}_j}-T_H$ has finite rank if and only if $H$ extends to a differentiable function on $\C\backslash\{0\}$ that solves the system of differential equations:
\begin{align*}
\frac{\partial H }{\partial z} &=\sum_{j=1}^NP'_j(z)\,\overline{Q}_j(\frac{1}{\bar{z}})\\
\frac{\partial H}{\partial \bar {z}}&=\sum_{j=1}^NP_j(\frac{1}{\bar{z}})\,\overline{Q'_j}(z)
\end{align*}
subject to the constraint $H(1)=\sum_{j=1}^NP_j(1)\,\overline{Q}_j(1)$.
\end{theoremB}

The paper is organized as follows. In Section 2, we recall the Mellin transform and use it to describe Toeplitz operators with quasihomogeneous symbols. We define the convolution $\diamond$ and prove a general result that implies Theorem A. In Section 3, we consider holomorphic polynomial symbols and show that our problem is equivalent to the existence of integrable solutions to a system of first order partial differential equations, providing a proof of Theorem B. Lastly, Section 4 provides a method to construct examples of polynomials $P_1,\ldots, P_N$ and $Q_1,\ldots, Q_N$ for which $\sum_{j=1}^{N}T_{P_j}T_{\overline{Q}_j}$ has finite rank. One of our examples actually recovers \cite[Theorem~4.4]{ding2017theorem}.

\section{Finite rank perturbations of Toeplitz products}

We first recall the definition of the Mellin transform, which plays an important role in our approach.
 
\begin{definition}
\label{Defn:MellinTransform}
The Mellin transform of a measurable function $\varphi$ on $(0,1)$ is defined as
\[\what{\varphi}(\zeta)=\int_0^1\varphi(r)\,r^{\zeta-1}dr\] for all complex numbers $\zeta$ for which the integral exists. 
\end{definition}

Note that if $\what{\varphi}(\alpha)$ exists for some complex number $\alpha$, then $\what{\varphi}(\zeta)$ is holomorphic on the open right half-plane $\{\zeta\in\C:\Re(\zeta)>\Re(\alpha)\}$ and is continuous on the closed right half-plane $\{\zeta\in\C:\Re(\zeta)\geq\Re(\alpha)\}$. If $\varphi\in L^1((0,1),rdr)$, then $\alpha$ can be chosen to be $2$. It can be seen from the definition that $\what{r^{m}\varphi}(\zeta)=\what{\varphi}(\zeta+m)$ whenever one of the sides is defined. This useful identity enables us to shorten several equations. See \cite{ButzerJFAA1997} for an exposition on the Mellin transform and its applications.

Recall that the Bergman space is a reproducing kernel Hilbert space with kernel $K(z,w)=(1-\bar{w}z)^{-2}$. It follows that $T_{f}$ is an integral operator. Indeed, for $h\in A^2$, we have
\begin{align}
\label{Eqn:int_formula_TO}
T_f(h)(z) & = \int_{\D}\frac{f(w)h(w)}{(1-\bar{w}z)^2}\,dA(w)\ \text{ for } z\in\disk.
\end{align}
When $f$ is unbounded but belongs to $L^1(\disk,dA)$, the integral on the right-hand side of \eqref{Eqn:int_formula_TO} may still exist, for example, when $h$ is bounded. One then uses \eqref{Eqn:int_formula_TO} to define the Toeplitz operator $T_{f}$ for $L^1$-symbol $f$. Such an operator may not be bounded on $A^2$ in general. However, if $f$ vanishes outside a disk $\{z: |z|\leq r\}$ for some $0<r<1$, then formula \eqref{Eqn:int_formula_TO} shows that for $h\in A^2$ and $z\in\D$,
\[|T_{f}(h)(z)| \leq \frac{\|f\|_{1}}{(1-r)^2}\big(\sup_{|w|\leq r}|h(w)|\big).\]
It then follows that $T_{f}$ is not only bounded but also compact. Consequently, if $f$ belongs to $L^1(\disk,dA)$ and is bounded on an annulus $\{z: r<|z|<1\}$ for some $0<r<1$, then $T_{f}$ is bounded. 

Mellin transform was used in the study Toeplitz operators by {\v{C}}u{\v{c}}kovi{\'c} and Rao \cite{CuckovicJFA1998} and has been then employed by many other researchers. Here we use Mellin transform to describe the action of Toeplitz operators with quasihomogeneous symbols on monomials. Let $m$ be an integer and $\varphi\in L^1((0,1),rdr)$. A direct calculation using \eqref{Eqn:int_formula_TO} and the series expansion of the kernel shows
\begin{equation}
\label{Eqn:quasihom_TO}
T_{\varphi(r) e^{im\theta}}(z^k)=
\begin{cases}
0, & k+m<0, \\
(2k+2m+2)\what{\varphi}(2k+m+2)\,z^{m+k}, & k+m\geq 0.
\end{cases}
\end{equation} 
Since $\{\sqrt{k+1}\,z^{k}: k=0,1,\ldots\}$ is an orthonormal basis for $A^2$, formula \eqref{Eqn:quasihom_TO} shows that $T_{\varphi(r)e^{im\theta}}$ is a weighted forward shift if $m\geq 0$ and a weighted backward shift if $m<0$.

Formula \eqref{Eqn:quasihom_TO} shows that for any functions $\psi$ and $\phi$ in $L^1((0,1),rdr)$, the Toeplitz product $T_{\phi(r) e^{im\theta}}T_{\psi(r) e^{in\theta}}$ is defined on holomorphic polynomials. Theorem 6.1 in \cite{louhichi2006products} provides criteria for such a product to equal another Toeplitz operator. The following result can be obtained by a similar argument which makes use of formula \eqref{Eqn:quasihom_TO}.

\begin{proposition}
\label{P:finiterank_LSZ}
Let $\phi,\psi$ and $\omega$ be functions in $L^1((0,1), rdr)$ and  $m,n$ be integers. Then following statements are equivalent.
\begin{enumerate}[(a.)]
\item The operator $T_{\phi(r) e^{im\theta}} T_{\psi(r) e^{in\theta}}-T_{\omega(r)e^{i(m+n)\theta}}$ has finite rank.
\item The equation $\widehat{\omega}(\zeta)=(\zeta-m+n)\widehat{\phi}(\zeta+n)\widehat{\psi}(\zeta-m)$ holds on some right half-plane in $\mathbb{C}$.
\end{enumerate} 
\end{proposition}

\begin{proof}
We sketch here a proof. Define $S=T_{\phi(r) e^{im\theta}} T_{\psi(r) e^{in\theta}}-T_{\omega(r)e^{i(m+n)\theta}}$. For all non-negative integers $k$ with $k+n\geq 0$ and $k+n+m\geq 0$, using formula \eqref{Eqn:quasihom_TO}, we see that $S(z^k)=(2k+2m+2n+2)\lambda_k\, z^{k+n+m}$, where
\[\lambda_k =(2k+2n+2)\what{\phi}(2k+2n+m+2)\what{\psi}(2k+n+2)-\what{\omega}(2k+m+n+2).\]
It follows that $S$ has finite rank if and only if $\lambda_k=0 $ for all but finitely many $k$. That is,
\begin{equation}
\label{Eqn:finite_rank_condition}
\what{\phi}(\zeta+n)\,\what{\psi}(\zeta-m) = \frac{1}{\zeta-m+n}\what{\omega}(\zeta)
\end{equation}
for all but finitely many $\zeta$ of the form $\zeta = 2k+m+n+2$. Since both sides are holomorphic and bounded for sufficiently large $\Re(\zeta)$, condition \eqref{Eqn:finite_rank_condition} is equivalent to (b) (see, for example, \cite[p.~102]{RemmertSpringer1998}).
\end{proof}

\begin{remark}
The operator $S=T_{\phi(r) e^{im\theta}} T_{\psi(r) e^{in\theta}}-T_{\omega(r)e^{i(m+n)\theta}}$ is, a priori, defined only on the space of holomorphic polynomials.
From the proof of Proposition \ref{P:finiterank_LSZ} we see that whenever $S$ has finite rank, it extends to a bounded operator on the whole space $A^2$.
\end{remark}

\begin{definition} We call a measurable function $\varphi$ defined on $(0,1)$ \textit{nearly integrable} if $\varphi$ belongs to $L^1((0,1),r^{\alpha}dr)$ for some positive number $\alpha$. Similarly, a measurable function $f$ defined on $\D$ is nearly integrable if $f$ belongs to $L^1(\D,|z|^{\alpha}dA)$ for some positive number $\alpha$. 
\end{definition}

Motivated by Proposition \ref{P:finiterank_LSZ}, we define a new convolution of functions that are linear combinations of quasihomogeneous functions.

\begin{definition}
\label{D:new_convolution} 
(a) Let $f_m(re^{i\theta})=\phi(r)e^{im\theta}$ and $g_n(re^{i\theta})=\psi(r)e^{in\theta}$ be nearly integrable quasihomogeneous functions of degrees $m$ and $n$ respectively. If there exists a nearly integrable function $\omega$ on $(0,1)$ such that
\begin{equation}
\label{Eqn:new_convolution}
\ds \widehat{\omega}(\zeta)=(\zeta-m+n)\,\widehat{\phi}(\zeta+n)\,\widehat{\psi}(\zeta-m)
\end{equation}
on a right half-plane in $\C$, then we define $f_m\diamond g_n$ as
\[\big(f_m\diamond g_n\big)(re^{i\theta})=\omega(r)e^{i(m+n)\theta}.\]
Note that for any complex constants $c$ and $d$, we have $(c f_m)\diamond (d g_m) = (cd)(f_m\diamond g_m)$.
(b) Let $\ds f=\sum_{m}f_m$ and $\ds g=\sum_{n}g_{n}$ be finite sums of nearly integrable quasihomogeneous functions. If for each pair $(m,n)$, the function $f_m\diamond g_{n}$ exists as in (a), then we define
\[f\diamond g=\sum_{m,n}f_m\diamond g_{n}.\]
\end{definition}

We first discuss two simple but important properties of the operation $\diamond$.

\begin{lemma}
\label{L:properties_diamond}
(a) The operation $\diamond$ is bilinear.

(b) Let $f$ and $g$ be as in Definition \ref{D:new_convolution} so that $f\diamond g$ is defined. Then for any integer $M\geq 0$, we have
$f\diamond (g\cdot z^M)=(f\diamond g)z^M$ and $(\bar{z}^M\cdot f)\diamond g=\bar{z}^{M}(f\diamond g)$. 
\end{lemma}

\begin{proof}
The bilinearity of $\diamond$ follows directly from the definition.

To prove (b), it suffices to consider the case $f$ and $g$ are quasihomogeneous functions. Write $f(re^{i\theta})=\phi(r)e^{im\theta}$, $g(re^{i\theta})=\psi(r)e^{in\theta}$ and $f\diamond g= \omega(r)e^{i(m+n)\theta}$, where $\omega$ satisfies $\ds \what{\omega}(\zeta)=(\zeta-m+n)\what{\phi}(\zeta+n)\what{\psi}(\zeta-m)$ on a right half-plane. Then
\begin{align*}
\what{r^M\omega}(\zeta)&=\what{\omega}(\zeta+M) =(\zeta+M-m+n)\what{\phi}(\zeta+M+n)\,\what{\psi}(\zeta+M-m)\\
&=(\zeta-m+(M+n))\what{\phi}(\zeta+(M+n))\,\what{r^M\psi}(\zeta-m).
\end{align*}
It follows that $\ds f\diamond (g\cdot z^M)(re^{i\theta})=r^M\omega(r)e^{i(m+M+n)\theta}=(f\diamond g)z^M$. The other identity is proved in a similar fashion.
\end{proof}

We now calculate the convolution $f\diamond g$ for a certain special class of functions.

\begin{example}
\label{Ex:convolution_monomial_type}
Let $f(re^{i\theta})=r^pe^{im\theta}$ and $g(re^{i\theta})=r^qe^{in\theta}$ for real numbers $p, q$ and integers $m,n$. Put $\ell=p+n-q+m$.  Then for $\Re(\zeta)$ sufficiently large, we have
\begin{align*}
(\zeta-m+n)\,\what{r^{p}}(\zeta+n)\,\what{r^{q}}(\zeta-m) & =\frac{\zeta-m+n}{(\zeta+p+n)(\zeta+q-m)}\\
& = \begin{cases}
\frac{n-q}{\ell}\cdot\frac{1}{\zeta+q-m} + \frac{m+p}{\ell}\cdot\frac{1}{\zeta+p+n}, & \ell\neq 0\\
\\
\frac{1}{\zeta+q-m} + \frac{n-q}{(\zeta+q-m)^2}, & \ell=0
\end{cases}\\ \\
& =\what{\omega}(\zeta),
\end{align*}
where \[\omega(r)=\begin{cases}
\frac{n-q}{l}r^{q-m}+\frac{m+p}{l}r^{p+n}, & \ell\neq 0\\
\\
r^{q-m}-(n-q)\log(r)r^{q-m}, & \ell=0.
\end{cases}\]
We have used the following formulas for the Mellin transform:
\[\what{r^{\alpha}}(\zeta)=\frac{1}{\zeta+\alpha},\qquad \what{r^{\alpha}\log(r)} = \frac{-1}{(\zeta+\alpha)^2}\] for $\Re(\zeta+\alpha)>0$. 
Consequently, we obtain
\begin{equation}
\label{Eqn:convolution_formula_monomials}
(f\diamond g)(re^{i\theta})=
\begin{cases}
\Big(\frac{n-q}{l}r^{q-m}+\frac{m+p}{l}r^{p+n}\Big)e^{i(m+n)\theta}, & \ell\neq 0\\
\\
r^{q-m}\Big(1-(n-q)\log(r)\Big)e^{i(m+n)\theta}, & \ell=0.\\
\end{cases}
\end{equation}
\end{example}

\begin{example}
\label{Ex:convolution_monomial}
Let $f(z)=z^{m}$ and $g(z)=\bar{z}^k$ for some integers $m$ and $k$. Applying the previous example with $p=m$ and $q=-n=k$ (and hence $\ell=2(m-k)$), we have
\begin{align*}
(z^m\diamond \bar{z}^{k}) 
& = \begin{cases}
\frac{n}{m-k}r^{k-m}e^{i(m-k)\theta} + \frac{m}{m-k}r^{m-k}e^{i(m-k)\theta}, & m-k\neq 0\\ \\
1+2k\log(r), & m-k=0
\end{cases}\\ \\
& = \begin{cases}
\frac{-k}{m-k}\bar{z}^{k-m} + \frac{m}{m-k}z^{m-k}, & m\neq k\\ \\
1+k\log|z|^2, & m=k.
\end{cases}
\end{align*}
Using the above formula, we compute the convolution of two Laurent polynomials in $z$ and $\bar{z}$. Let
\[f(z) = \sum_{m,n}a_{mn}z^{m}\bar{z}^{n}\qquad\text{and}\qquad g(z) = \sum_{\ell,k}b_{\ell k}z^{\ell}\bar{z}^{k}\]
be two Laurent polynomials (that is, the above sums are over finite sets of integers). Then the convolution $f\diamond g$ exists and is given by
\begin{align*}
f\diamond g & = \sum_{m,n,\ell,k}a_{mn}b_{\ell k}(z^{m}\bar{z}^{n})\diamond (z^{\ell}\bar{z}^{k})\\
& = \sum_{m,n,\ell,k}a_{mn}b_{\ell k}\,\bar{z}^{n}\big(z^m\diamond \bar{z}^{k}\big){z}^{\ell}\quad\text{(by Lemma \ref{L:properties_diamond}(b))}\\
& = \sum_{\substack{m,n,\ell,k\\ m\neq k}}a_{mn}b_{\ell k}\,\bar{z}^{n}\Big(\frac{-k}{m-k}\bar{z}^{k-m}+\frac{m}{m-k}z^{m-k}\Big){z}^{\ell} \\
& \qquad\qquad + \sum_{m,n,\ell}a_{mn}b_{\ell m}\,\bar{z}^{n}\big(1+m\log|z|^2\big)z^{\ell}\\
& = \sum_{\substack{m,n,\ell,k\\ m\neq k}}a_{mn}b_{\ell k}\Big\{\frac{-k}{m-k}z^{\ell}\bar{z}^{n+k-m}+\frac{m}{m-k}z^{m+\ell-k}\bar{z}^{n}\Big\}\\
& \qquad\qquad + \sum_{m,n,\ell}a_{mn}b_{\ell m}\,\bar{z}^{n}\big(1+m\log|z|^2\big)z^{\ell}.
\end{align*}
\end{example}

It is a natural problem to investigate the existence of $f\diamond g$ for general quasihomogeneous functions $f$ and $g$. From Definition \ref{D:new_convolution}, the existence of $\varphi(r)e^{im\theta}\diamond \psi(r)e^{in\theta}$ depends on the existence of a nearly integrable function $\omega$ on the interval $(0,1)$ that satisfies condition \eqref{Eqn:new_convolution}. This is equivalent to the requirement that the inverse Mellin transform of $F(\zeta)=(\zeta-m+n)\widehat{\varphi}(\zeta+n)\widehat{\psi}(\zeta-m)$ exists and is supported on $(0,1)$. Obtaining necessary and sufficient conditions on such functions $\varphi$ and $\psi$ seems to be a delicate problem. We offer here a sufficient condition, which covers a wide range of cases.

\begin{proposition}
\label{P:existence_new_convolution} Let $\phi,\psi$ be nearly integrable functions on $(0,1)$ and $m,n$ be integers. Assume that $\phi$ (or $\psi$) is differentiable on $(0,1]$ and $\phi'$ (respectively, $\psi'$) is nearly integrable. Then the convolution $\phi(r)e^{im\theta}\diamond \psi(r)e^{in\theta}$ exists. \\
As a consequence, if
\[f(re^{i\theta}) = \sum_{j=1}^{d}\phi_j(r)e^{im_j\theta},\qquad g(re^{i\theta}) = \sum_{j=1}^{d}\psi_j(r)e^{in_j\theta}\] are finite sums of nearly integrable quasihomogeneous functions such that for each $j$, the function $\phi_j$ (or $\psi_j$) is differentiable on $(0,1]$ and $\phi'$ (respectively, $\psi'$) is nearly integrable, then $f\diamond g$ exists.
\end{proposition}

\begin{proof}
From the assumption on $\phi$, we see that for all complex numbers $\zeta$ with sufficiently large real part, the function 
\[h(r) = \begin{cases}
r^{\zeta+n}\phi(r), & \quad 0<r\leq 1\\
0, & \quad r = 0
\end{cases}\]
is differentiable on $[0,1]$ and $h'$ is integrable. Writing $h(r)=r^{\zeta-m+n}\cdot (r^m\phi(r))$, we compute the derivative
\[h'(r) = (\zeta-m+n)\phi(r)r^{\zeta+n-1} + r^{\zeta-m-1}\varphi(r),\]
where $\varphi = r^{n+1}d(r^m\phi(r))/dr$. Then
\begin{align*}
\phi(1) = h(1)-h(0) = \int_{0}^{1} h'(r)dr = (\zeta-m+n)\widehat{\phi}(\zeta+n) +  \widehat{\varphi}(\zeta-m).
\end{align*}
It follows that
\begin{align*}
(\zeta-m+n)\widehat{\phi}(\zeta+n)\widehat{\psi}(\zeta-m) & = \phi(1)\widehat{\psi}(\zeta-m) - \widehat{\varphi}(\zeta-m)\widehat{\psi}(\zeta-m).
\end{align*}
As consequence, equation \eqref{Eqn:new_convolution} admits a solution $\omega$ given by
\[\omega(r) = r^{-m}\Big(\phi(1)\psi(r)-(\varphi*\psi)(r)\Big),\quad 0<r<1,\]
where $\varphi*\psi$ denotes the Mellin convolution defined as
\[(\varphi*\psi)(r) = \int_{r}^{1}\varphi(t)\psi(\frac{r}{t})\frac{dt}{t}.\]
We have used the well-known fact that $\widehat{\varphi*\psi}=\widehat{\varphi}\cdot\widehat{\psi}$.
\end{proof}

Using Proposition \ref{P:finiterank_LSZ}, we now show that the convolution $\diamond$ is indeed important for our study of finite rank perturbations of Toeplitz products. We shall consider Toeplitz operators whose symbols belong to the linear span of integrable quasihomogeneous functions. Let us denote this class of symbols by \SQH (span of quasihomogeneous functions), that is,
\[\SQH = \Big\{\sum_{m\in E}f_m: \ \text{ finite } E\subset\mathbb{Z}\ \text{ and } f_m(re^{i\theta})=\phi_m(r)e^{im\theta}\in L^1(\D,dA), \forall m\in E\Big\}.\]

\begin{theorem}
\label{T:new_convolution}
(a) Let $f,g$ belong to \SQH such that $f\diamond g$ is defined. Then there exists an integer $M\geq 0$ such that $T_{f}T_{g\cdot z^M}-T_{(f\diamond g)z^M}$ is a finite rank operator.

(b) Let $F_1,\ldots, F_N$ and $G_1,\ldots, G_N$ belong to \SQH  such that  $F_j\diamond G_j$ is defined for all $1\leq j\leq N$. Let $H\in L^1(\D,dA)$. Then $\sum_{j=1}^{N} T_{F_j}T_{G_j} - T_H$ has finite rank if and only if $\sum_{j=1}^{N}F_j\diamond G_j$ belongs to $L^1(\D,dA)$ and $H=\sum_{j=1}^{N}F_j\diamond G_j$.

In the case $F_j$'s and $G_j$'s are polynomials in $z$ and $\bar{z}$, we obtain Theorem A.
\end{theorem}

\begin{proof}
(a) Write ${f=\sum_{m\in E} f_m}$ and $g=\sum_{n\in F} g_{n}$, where $f_m(re^{im\theta})=\phi_m(r)e^{im\theta}$ and $\ds g_{n}(re^{i\theta})=\psi_n(r)e^{in\theta}$ belong to $L^1(\D,dA)$. By the definition of $\diamond$, for each pair $(m,n)\in E\times F$, the function $f_m\diamond g_n$ is defined and is nearly integrable (but may not be integrable). There then exists an integer $M\geq 0$ such that $(f_m\diamond g_n)z^M$ belongs to  $L^1(\D,dA)$ for all $(m,n)\in E\times F$. Using Lemma \ref{L:properties_diamond}, we have
\begin{align*}
T_{f}T_{g\cdot z^M} - T_{(f\diamond g)z^M} & = T_{f}T_{g\cdot z^M} - T_{f\diamond (g\cdot z^M)} = \sum_{m,n}\big(T_{f_m}T_{g_n\cdot z^M} - T_{f_m\diamond (g_n \cdot z^M)}\big).
\end{align*}
By Proposition \ref{P:finiterank_LSZ}, the operator $T_{f_m}T_{g_n\cdot z^M} - T_{f_m\diamond (g_n\cdot z^M)}$ has finite rank. As a consequence, $T_{f}T_{g\cdot z^M} - T_{(f\diamond g)z^M}$, being a finite sum of finite-rank operators, has finite rank as well.

(b) Put $h=\sum_{j=1}^{N}F_j\diamond G_j$. Using (a), we can find an integer $M\geq 0$ such that the operator $T=\sum_{j=1}^{N} T_{F_j}T_{G_j\cdot z^M} - T_{h\cdot z^M}$ has finite rank. 

If $h$ belongs to $L^1(\D,dA)$, then $T_{h}$ is defined on the space of holomorphic polynomials and $T_{h\cdot z^M}=T_{h}T_{z^M}$. It follows that
\[T = \Big(\sum_{j=1}^{N} T_{F_j}T_{G_j} - T_{h}\Big)T_{z^M}.\]
Since $T$ has finite rank and the range of $T_{z^M}$ has finite codimension, we conclude that $\sum_{j=1}^{N} T_{F_j}T_{G_j} - T_{h}$ has finite rank.
 
Now suppose that $\ds S=\sum_{j=1}^{N}T_{F_j}T_{G_j}-T_{H}$
has finite rank. Then
\[ST_{z^M} = \sum_{j=1}^{N}T_{F_j}T_{G_j\cdot z^M} - T_{H\cdot z^M}\] also has finite rank. It follows that
\[T_{(h-H)\cdot z^M} = T_{h\cdot z^M} - T_{H\cdot z^M} = T - ST_{z^M}\] is a finite rank Toeplitz operator. By Luecking's Theorem \cite{luecking2008finite}, we conclude that $(h-H)\cdot z^M=0$ a.e., which implies $h=H$ a.e. on the unit disk. That is, $\sum_{j=1}^{N}F_j\diamond G_j = H$, which completes the proof of the theorem.
\end{proof}

\begin {corollary}
\label{C:complete_classification_type_I}
Let $f(re^{i\theta})=r^pe^{im\theta}$ and $g(re^{i\theta})=r^qe^{in\theta}$ belong to $L^1(\D,dA)$, where $m,n$ are integers and $p>-2$ and $q>-2$ are real numbers. Then the following statements are equivalent.

\begin{enumerate}[(a.)]
\item There exists $h\in L^1(\mathbb{D}, dA)$ such that  $T_{f}T_{g}-T_h$  is of finite rank.
\item One of the following conditions holds true.
 \begin{enumerate}[(1.)]
\item  $ q-m+2 > 0$ and $p+n+2>0$;
\item  $m+p=0$ and $p+n+2>0$;
\item  $n-q=0$ and $q-m+2>0$.
\end{enumerate}
\end{enumerate}
If one of the statements, hence both hold true, then $h=f\diamond g$. 
\end{corollary}

\begin{proof}
Theorem \ref{T:new_convolution} shows that (a) is equivalent to the condition that $f\diamond g$ belongs to $L^1(\D,dA)$. Since $r^{\alpha}$ and $r^{\alpha}\log(r)$ belong to $L^1((0,1),rdr)$ if and only if $\alpha+2>0$, formula \eqref{Eqn:convolution_formula_monomials} shows that  $f\diamond g$ belongs to $L^1(\D,dA)$ if and only if (b) holds.
\end{proof}

\begin{remark}
In Example \eqref{ExampleA}, we have $p=m=2$, $q=3$ and $n=-3$. On the other hand, in Example \eqref{ExampleB}, we have $p=m=2$, $q=4$ and $n=-2$.
\end{remark}

\section{Holomorphic polynomial symbols}

\v{C}u\v{c}kovi\'{c} \cite{cuckovicIEOT2007} obtained necessary and sufficient conditions for which $T_{f}T_{g}-T_{h^n}$ is a finite rank operator under the assumption that $f, g$ and $h$ are bounded harmonic functions. In a recent paper, Ding, Qin and Zheng studied a more general problem by replacing $h^n$ by an arbitrary bounded $C^2$-function whose invariant Laplacian is integrable. They showed \cite[Theorem 1.6]{ding2017theorem} that in the case the rank $r\leq 1$, either $\bar{f}$ or $g$ must be holomorphic and as a consequence, $r=0$. On the other hand, \cite[Theorem 1.8]{ding2017theorem} provided a complete characterization in the case of higher ranks.

Note that for $f=f_1+\bar{f}_2$ and $g=g_1+\bar{g}_2$, where the functions $f_1, f_2, g_1, g_2$ are bounded holomorphic, the Toeplitz product $T_{f}T_{g}$ is a finite rank perturbation of another Toeplitz operator if and only if the same holds true for the product $T_{f_1}T_{\bar{g}_2}$. In this section, we shall investigate the question: suppose $P_1,\ldots, P_N$ and $Q_1,\ldots, Q_N$ are holomorphic polynomials and $H$ belongs to $L^1(\D,dA)$, under which conditions does the operator $\sum_{j=1}^{N}T_{P_j}T_{\overline{Q}_j}-T_{H}$ have finite rank? Even though we have such a restriction on the symbols, we do allow $H$ to be non-differentiable and even unbounded. We also consider finite sums of several Toeplitz products. The approach in \cite{ding2017theorem}, which makes use of Axler-Chang-Sarason's result \cite{AxlerChangSarasonIEOT1978} on finite rank semicommutators of Toeplitz operators on the Hardy space, does not seem to work in our settings. It is interesting and surprising that our results in this section reveal a connection between the theory of Toeplitz operators and solutions to a certain system of first order differential equations.

\begin{definition}
Recall that a \textit{Laurent polynomial} is a finite linear combination of $\{z^{m}: m\in\mathbb{Z}\}$ with complex coefficients. Such a polynomial is holomorphic on the punctured complex plane $\C\backslash\{0\}$.
\end{definition}

\begin{proposition}
\label{P:diamond_solution_DE}
Let $P_j$ and $Q_j$ be Laurent polynomials for $j=1,\ldots, N$. Then the function $G=\sum_{j=1}^{N}P_j\diamond\overline{Q}_j$ is the unique solution to the following system of differential equations on $\C\backslash\{0\}$:
\begin{align*}
\frac{\partial G }{\partial z} &=\sum_{j=1}^NP'_j(z)\,\overline{Q}_j(\frac{1}{\bar{z}})\\
\frac{\partial G}{\partial \bar {z}}&=\sum_{j=1}^NP_j(\frac{1}{\bar{z}})\,\overline{Q'_{j}}(z)
\end{align*}
subject to the constraint $G(1)=\sum_{j=1}^NP_j(1)\,\overline{Q}_j(1)$.
\end{proposition}

\begin{proof}
For integers $m$ and $k$, define $f_m(z)=z^m$ and $g_k(z)=z^{k}$. Example \ref{Ex:convolution_monomial} gives
\[(f_m\diamond \bar{g}_k)(z) = \begin{cases}
\frac{-k}{m-k}\bar{z}^{k-m} + \frac{m}{m-k}z^{m-k}, & m\neq k\\ \\
1+k\log|z|^2, & m=k.
\end{cases}\]
We see that $(f_m\diamond \bar{g}_k)(1)=1 = f_m(1)\bar{g}_k(1)$ and in both cases,
\[\frac{\partial (f_m\diamond \bar{g}_k) }{\partial z} =
mz^{m-k-1} = f'_m(z)\,g_k\big(\frac{1}{z}\big),\]
and
\[\frac{\partial (f_m\diamond \bar{g}_k) }{\partial \bar {z}} =
k\bar{z}^{k-m-1} = f_m\big(\frac{1}{\bar{z}}\big)\,\overline{g'_k}(z).\]
Now let $P(z)=\sum_{m=-M}^{M} a_mz^m$ and $Q(z)=\sum_{k=-K}^K b_kz^k$ be two Laurent polynomials. Using the linearity of $\diamond$, we have 
\[P\diamond \overline{Q}=\sum_{m,k}a_m\bar{b}_kf_m\diamond \bar{g}_k.\]
It then follows that
\[ \frac{\partial (P\diamond \overline{Q}) }{\partial \bar {z}}=\sum_{m,k}a_m\bar{b}_k\frac{\partial (f_m\diamond \bar{g}_k) }{\partial \bar {z}}=\sum_{m,k}a_m\bar{b}_kf_m\big(\frac{1}{\bar{z}}\big)\,\overline{g'_k}(z)=P(\frac{1}{\bar{z}})\,\overline{Q'}(z).\]
Similarly,
\[ \frac{\partial (P\diamond \overline{Q}) }{\partial z}=P'(z)\overline{Q}(\frac{1}{\bar{z}}).\]
Also,
 \[(P\diamond \overline{Q})(1)=\sum_{m,k}a_m\bar{b}_k\,(f_m\diamond \bar{g}_k)(1)=\sum_{m,k} a_m\bar{b}_k=P(1)\,\overline{Q}(1).\]
The proof is completed for the case $N=1$. The general case for arbitrary positive integers $N$ follows from linearity. 

Finally, uniqueness of the solution is standard and it follows from the fact that $W=0$ is the only solution to the system
\[\frac{\partial W}{\partial z} = \frac{\partial W}{\partial\bar{z}} = 0\]
with initial condition $W(1)=0$.
\end{proof}

We are now ready for a proof of Theorem B, which we restate below for the reader's convenience.

\begin{theorem}
\label{T:holomorphic_polynomials}
Let  $P_j$ and $Q_j$ be Laurent polynomials for $j=1,\ldots, N$ and let $H$ belong to  $L^1(\D,dA)$. Then $\ds\sum_{j=1}^{N}T_{P_j}T_{\overline{Q}_j}-T_H$ has finite rank if and only if $H$ extends to a differentiable function on $\C\backslash\{0\}$ that solves the system of differential equations:
\begin{align*}
\frac{\partial H }{\partial z} &=\sum_{j=1}^NP'_j(z)\,\overline{Q}_j(\frac{1}{\bar{z}})\\
\frac{\partial H}{\partial \bar {z}}&=\sum_{j=1}^NP_j(\frac{1}{\bar{z}})\,\overline{Q'_j}(z)
\end{align*}
subject to the constraint $H(1)=\sum_{j=1}^NP_j(1)\,\overline{Q}_j(1)$.
\end{theorem}

\begin{proof}
It follows from Theorem \ref{T:new_convolution} that the operator $\ds\sum_{j=1}^{N}T_{P_j}T_{\overline{Q}_j}-T_H$ has finite rank if and only if $H=\sum_{j=1}^{N}P_j\diamond \overline{Q}_j$, which, by Proposition \ref{P:diamond_solution_DE}, is equivalent to the requirement that $H$ solves the given system of partial differential equations with the constraint.
\end{proof}

In the case $N=1$, we are able to describe more explicitly which polynomials satisfy the conditions in Theorem \ref{T:holomorphic_polynomials}.

\begin{corollary}
\label{C:single_product_polynomials}
Let  $P$ and $Q$ be two holomorphic polynomials in $z$. Then there exists $H\in L^1(\D,dA)$ such that $T_PT_{\bar{Q}}-T_H$ has finite rank if and only if one of the following conditions is satisfied:
\begin{enumerate}[(1)]
\item $P$ is a constant or $Q$ is a constant.
\item $P(z)=a_0 + a_{M-2}z^{M-2}+a_{M-1}z^{M-1}+a_Mz^{M}$ and  $Q(z)=b_0 + b_{M-1}z^{M-1}$ with $a_Mb_{M-1}\neq 0$.
\item $P(z) = a_0 + a_{M-1}z^{M-1}+a_{M}z^M$ and $Q(z)=b_0+b_{M-1}z^{M-1}+b_{M}z^M$ with $a_Mb_M\neq 0$.
\item $P(z)=a_0 +a_Mz^{M}$ and $Q(z)=b_0 + b_{M-1}z^{M-1}+b_Mz^{M}+b_{M+1}z^{M+1}$ with $a_Mb_{M+1}\neq 0$.
\end{enumerate}
\end{corollary}

\begin{proof}
We only need to investigate the existence of $L^1$-solution to the system of differential equations
\begin{align*}
\frac{\partial H }{\partial z} =P'(z)\,\overline{Q}(\frac{1}{\bar{z}}),\qquad
\frac{\partial H }{\partial \bar {z}}=P(\frac{1}{\bar{z}})\,\overline{Q'}(z)
\end{align*}
with constraint $H(1)=P(1)\overline{Q}(1)$. It is clear that if either $P$ or $Q$ is a constant, then $H=P\bar{Q}$ is a solution.

Assume now that both $P$ and $Q$ are non-constant polynomials. Write $P(z)=a_0+\sum_{k=m}^Ma_kz^k$ and  $Q(z)=b_0+\sum_{j=n}^Nb_jz^j$, where $1\leq m\leq M$, $1\leq n\leq M$, $a_m, a_M$ and $b_n, b_N$ are non-zero complex numbers. Then 
\begin{align*}
P'(z)\,\overline{Q}(\frac{1}{\bar{z}}) &=\Big(\sum_{k=m}^Mka_kz^{k-1}\Big)\Big(\bar{b}_0+\sum_{j=n}^N\bar{b}_jz^{-j}\Big) = ma_m\bar{b}_{N}z^{m-1-N} + \text{higher powers},\\
P(\frac{1}{\bar{z}})\,\overline{Q'}(z)&=\Big(a_0+\sum_{k=m}^Ma_k\bar{z}^{-k}\Big)\Big(\sum_{j=n}^Nj\bar{b}_j\bar{z}^{j-1}\Big) = na_M\bar{b}_n\bar{z}^{n-1-M} +\text{higher powers}.
\end{align*}
We see that any solution $H$ must be of the form
\[H(z) = H_1(z) + \overline{H_2(z)} + \alpha\log|z|^2 + \beta,\] where $\alpha,\beta$ are complex constant, $H_1$ is a Laurent polynomial with lowest power $z^{m-N}$ and $H_2$ is a Laurent polynomial with lowest power $z^{n-M}$. In order for $H$ to belong to $L^1(\D,dA)$, it is necessary and sufficient that $m-N\geq -1$ and $n-M\geq -1$, which means $m\geq N-1$ and $n\geq M-1$. Since $1\leq m\leq M$ and $1\leq n\leq N$, we conclude that
\[N-1\leq m\leq M\ \text{ and }\ M-1\leq n\leq N.\]
It follows that $-1\leq N-M\leq 1$ and hence there are three cases to consider.

If $N=M-1$, then $M-2\leq m\leq M$ and $n=M-1$. This gives condition (2).

If $N=M$, then $M-1\leq m, n \leq M$. This gives condition (3).

If $N=M+1$, then $m=M$ and $M-1\leq n\leq M+1$. This gives condition (4).
\end{proof}

Using Theorem \ref{T:holomorphic_polynomials} in the case $H=0$, we obtain necessary and sufficient conditions for a finite sum of products of two Toeplitz operators with rational symbols having poles outside of the closed unit disk to have finite rank. In the case of polynomial symbols, our conditions are equivalent to those obtained in \cite[Theorem 4.2]{ding2017theorem}. Their approach, which uses Berezin transform, is different from ours.

\begin{theorem}  
\label{T:generalized_DQZ}
Let  $F_j$ and $G_j$ be rational functions whose poles lie outside of the closed unit disk for $j=1,\ldots,N$. Then the operator $\sum_{j=1}^{N}T_{F_j}T_{\overline{G}_j}$ has finite rank if and only if for all but a discrete set of $z\in\C$,
\begin{align}
\sum_{j=1}^{N}F_j(z)\,\overline{G}_j(\frac{1}{\bar{z}}) & = 0, \label{Eqn:DQZ1}
\end{align}
and
\begin{align}
\sum_{j=1}^{N}F_j(z)\,\overline{G'_j}(\frac{1}{\bar{z}}) & =0.\label{Eqn:DQZ2}
\end{align}
\end{theorem}

\begin{proof}
Let $D$ denote the common denominator for all $F_j$'s and $G_j$'s. Then there are polynomials $P_1,\ldots, P_N$ and $Q_1,\ldots,Q_N$ such that $DF_j = P_j$ and $DG_j = Q_j$ for all $1\leq j\leq N$. Using the identities $T_{P_j}=T_{D}T_{F_j}$ and $T_{\overline{Q}_j} = T_{\overline{G}_j}T_{\overline{D}}$, we have
\begin{align*}
\sum_{j=1}^{N}T_{P_j}T_{\overline{Q}_j} & = T_{D}\Big(\sum_{j=1}^{N}T_{F_j}T_{\overline{G}_j}\Big)T_{\overline{D}}
\end{align*}

Since $T_{D}$ is injective and $T_{\overline{D}}=T_{D}^{*}$ has dense range, we see that $\sum_{j=1}^{N}T_{F_j}T_{\overline{G}_j}$ has finite rank if and only if 
$\sum_{j=1}^{N}T_{P_j}T_{\overline{Q}_j}$ has finite rank which, by Theorem \ref{T:holomorphic_polynomials} with $H=0$, is equivalent to
\begin{equation}
\label{Eqn:P_jQ_j}
\sum_{j=1}^{N}P_j(1)\overline{Q}_j(1)=0\text{ and }\sum_{j=1}^NP_j'(z)\,\overline{Q}_j(\frac{1}{\bar{z}}) = \sum_{j=1}^N P_j(z)\,\overline{Q_j'}(\frac{1}{\bar{z}}) =0\text{ for } z\neq 0.
\end{equation}
Define
\[U(z) = \sum_{j=1}^{N}P_j(z)\overline{Q}_j(\frac{1}{\bar{z}}),\qquad z\neq 0.\]
Then $U$ is holomorphic on $\C\backslash\{0\}$ and 
\[U'(z) = \sum_{j=1}^{N}P_j'(z)\,\overline{Q}_j(\frac{1}{\bar{z}}) -\frac{1}{z^2}\sum_{j=1}^{N}P_j(z)\,\overline{Q_j'}(\frac{1}{\bar{z}}).\]
It follows that condition \eqref{Eqn:P_jQ_j} is equivalent to 
\begin{align}
\label{Eqn:U}
\sum_{j=1}^NP_j'(z)\,\overline{Q}_j(\frac{1}{\bar{z}}) = 0\quad\text{ and }\quad U(z)=0\text{ for } z\neq 0.
\end{align}
Using the fact that $P_j=DF_j$ and $Q_j = DG_j$, where $D$ is nonzero, it is not difficult to see that \eqref{Eqn:U} is equivalent to \eqref{Eqn:DQZ1} and \eqref{Eqn:DQZ2}.
\end{proof}

\section{Examples}

In this section, we discuss a method to construct examples of polynomials $P_1,\ldots, P_N$ and $Q_1,\ldots, Q_N$ such that the operator $S=\sum_{j=1}^{N}T_{P_j}T_{\overline{Q}_j}$ has finite rank. Assume that $P_j$'s and $Q_j$'s are of degrees at most $d$. For $1\leq j\leq N$, write $Q_j(z) = \sum_{\ell=0}^{d}q_{\ell,j}z^{\ell}$ and $P_{j}(z)=\sum_{k=0}^{d}p_{k,j}z^{k}$. We then have
\begin{align*}
S = \sum_{j=1}^{N}\sum_{0\leq k,\ell\leq d}p_{k,j}\bar{q}_{\ell,j}T_{z^k}T_{\bar{z}^{\ell}} = \sum_{0\leq k,\ell\leq d}\big(\sum_{j=1}^{N}p_{k,j}\bar{q}_{\ell,j}\big)T_{z^k}T_{\bar{z}^{\ell}} = \sum_{0\leq k,\ell\leq d}c_{k,\ell}T_{z^k}T_{\bar{z}^{\ell}},
\end{align*}
where $c_{k,\ell} = \sum_{j=1}^{N}p_{k,j}\bar{q}_{\ell,j}$. Put $\bfP=(p_{k,j})$, $\bfQ=(q_{\ell,j})$ and $\bfS = (c_{k,\ell})$. Then $\bfS = \bfP\cdot\bfQ^{*}$, where $\bfQ^{*}$ is the conjugate transpose of $\bfQ$.
By Theorem \ref{T:generalized_DQZ}, the operator $S$ has finite rank if and only if 
\begin{align*}
0 & = \sum_{0\leq k, \ell\leq d}c_{k,\ell}\,z^{k}\cdot\Big(\frac{1}{z^{\ell}}\Big) = \sum_{m=-d}^{d}\Big(\sum_{k-\ell=m}c_{k,\ell}\Big)z^{m}
\end{align*}
and
\begin{align*}
0 & = \sum_{0\leq k, \ell\leq d}c_{k,\ell}\,z^{k}\cdot\Big(\frac{-\ell}{z^{\ell+1}}\Big) = -\sum_{m=-d}^{d}\Big(\sum_{k-\ell=m}\ell c_{k,\ell}\Big)z^{m-1},
\end{align*}
which are equivalent to
\begin{align}
\label{Eqn:S_finite_rank_conditions}
\sum_{k-\ell=m}{c}_{k,\ell} = \sum_{k-\ell=m}\ell\,{c}_{k,\ell}=0 \quad\text{ for all }\quad m\in\{-d,\ldots,d\}.
\end{align}

\begin{example}
\label{E:sum_three_products}
Let $\mathbf{S_d}$ be the $(d+1)\times(d+1)$ matrix whose entries are all zero except $c_{0,0}=d-1$, $c_{1,1}=-d$, $c_{d,d}=1$. Then the conditions \eqref{Eqn:S_finite_rank_conditions} are satisfied so the operator
\begin{align}
\label{Eqn:S_canonical_representation}
S_d = (d-1)T_{1}T_{1} + (-d)T_{z}T_{\bar{z}} + T_{z^{d}}T_{\bar{z}^d}
\end{align}
has finite rank. A direct calculation using the formula
\begin{equation}
T_{z^m}T_{\bar{z}^m}(z^k)  = \begin{cases} 0 & \text{ if } 0\leq k\leq m-1,\\
\frac{k-m+1}{k+1}z^{k} & \text{ if } k\geq m,
\end{cases}
\end{equation}
shows that for $d\geq 2$,
$\rank(S_d)=d-1$ and
\[S_d = \sum_{k=0}^{d-2}\big(d-k-1\big)z^{k}\otimes z^k.\]
\end{example}

\begin{example}
\label{E:Theorem4.4DQZ}
In this example, we recover \cite[Theorem~4.4]{ding2017theorem}. Suppose $\alpha, \beta, \gamma$ are real numbers such that $\beta\neq 0$ and $\alpha\neq\gamma$. Take $\bfP$ and $\bfQ$ to be $(d+1)\times 3$  matrices whose entries are all zero except the submatrices formed by the $0$th, $1$st and $d$th rows are given as
\begin{align*}
\begin{pmatrix}
-(d-1)\gamma & -(d-1)\alpha & (d-1)(\alpha-\gamma)\\
d\beta & d\beta  & 0\\
\beta\gamma & \beta\alpha & 0
\end{pmatrix}
\quad\text{and}\quad
\begin{pmatrix}
0 & 0 & -\beta\\
\alpha & -\gamma & 0\\
1 & -1 & -1
\end{pmatrix}.
\end{align*}
It is clear that $\bfP\cdot \bfQ^{*} = \beta(\gamma-\alpha)\mathbf{S_d}$, where $\mathbf{S_d}$ is given in Example \ref{E:sum_three_products}. As a consequence, the operator $T_{F_1}T_{\overline{G}_1} - T_{F_2}T_{\overline{G}_2} - T_{F_3}T_{\overline{G}_3} = \beta(\gamma-\alpha)S_{d}$ has rank $d-1$. Here,
\begin{align*}
F_1(z) & = -(d-1)\gamma + d\beta z + \beta\gamma z^d,\\
F_2(z) & = -(d-1)\alpha + d\beta z + \beta\alpha z^d,\\
F_3(z) & = (d-1)(\alpha-\gamma)\\
G_1(z) & = \alpha z + z^d,\\
G_2(z) & = \gamma z + z^d,\\
G_3(z) & = \beta + z^d.
\end{align*}
These polynomials are constructed using the columns of $\mathbf{P}$ and $\mathbf{Q}$. For $d\geq 3$, under appropriate conditions on the constants $\alpha, \beta$ and $\gamma$ (see \cite[Theorem~4.4]{ding2017theorem}), the functions $F_2$ and $G_3$ do not vanish on the closed unit disk. It then follows that the operator
\[T_{\frac{F_1}{F_2}}T_{\frac{\overline{G}_1}{\overline{G}_3}} - T_{\frac{F_3}{F_2}+\frac{\overline{G}_2}{\overline{G}_3}} = T_{\frac{1}{F_2}}\Big(T_{F_1}T_{\overline{G}_1} - T_{F_2}T_{\overline{G}_2} - T_{F_3}T_{\overline{G}_3}\Big)T_{\frac{1}{\overline{G}_3}}\]
has rank $d-1$. 
\end{example}

\begin{example}
\label{E:revisedDQZ} In this example, we exhibit a different collection of polynomials that satisfy the same properties as in the previous example. Let $\beta$ be any real number. Take $\bfP$ and $\bfQ$ to be $(d+1)\times 3$ matrices whose entries are all zero except the submatrices formed by the $0$th, $1$st and $d$th rows are given as
\begin{align*}
\begin{pmatrix}
0 & (d-1) & -(d-1)\beta\\
-d(\beta^2+1) & 0 & 0\\
0 & -\beta & -1 &
\end{pmatrix}
\quad\text{and}\quad
\begin{pmatrix}
0 & 1 & -{\beta}\\
1 & 0 & 0\\
0 & -\beta & -1
\end{pmatrix}.
\end{align*}
\end{example}
Since $\bfP\cdot\bfQ^{*}=(\beta^2+1)\bfS_{d}$, if we define
\begin{align*}
F_1(z) & = -d(\beta^2+1)z,\\
F_2(z) & = (d-1) - \beta z^d,\\
F_3(z) & = -(d-1)\beta - z^d,\\
G_1(z) & = z,\\
G_2(z) & = -1+\beta z^d,\\
G_3(z) & = \beta + z^d,
\end{align*}
then the operator
$T_{F_1}T_{\overline{G}_1} - T_{F_2}T_{\overline{G}_2} - T_{F_3}T_{\overline{G}_3} = (\beta^2+1)S_{d}$ has rank $d-1$. If $d\geq 3$ then for all $1<|\beta|<d-1$, the functions $F_2$ and $G_3$ do not vanish in the closed unit disk so the same argument as in Example \ref{E:Theorem4.4DQZ} shows that the operator $T_{\frac{F_1}{F_2}}T_{\frac{\overline{G}_1}{\overline{G}_3}} - T_{\frac{F_3}{F_2}+\frac{\overline{G}_2}{\overline{G}_3}}$ has rank $d-1$.

\subsection*{Acknowledgements} The authors wish to thank the referee for helpful suggestions which improved the presentation of the paper.

\bibliographystyle{amsplain}

\providecommand{\bysame}{\leavevmode\hbox to3em{\hrulefill}\thinspace}
\providecommand{\MR}{\relax\ifhmode\unskip\space\fi MR }
\providecommand{\MRhref}[2]{%
  \href{http://www.ams.org/mathscinet-getitem?mr=#1}{#2}
}
\providecommand{\href}[2]{#2}

\end{document}